\documentclass{amsart}

\usepackage{amsmath, amssymb}
\usepackage{amsrefs}

\newtheorem{theorem}{Theorem}[section]
\newtheorem{lemma}[theorem]{Lemma}

\newtheorem{definition}[theorem]{Definition}

\newtheorem{remark}[theorem]{Remark}

\newcommand{\Out}{\mathrm{Out}}
\newcommand{\id}{\mathrm{id}}

\title{A criterion for Kolchin subgroups of $\Out(F_r)$}
\author{Edgar A. Bering IV}

\begin{document}

\maketitle

\begin{abstract}
This article provides a decidable criterion for when a subgroup of $\Out(F_r)$ generated
by two Dehn twists consists entirely of polynomially growing elements,
answering an earlier question of the author.
\end{abstract}
\section{Introduction}

Outer automorphims of a free group divide into two categories, polynomially
growing and exponentially growing, according to the behaviour of word lengths
under iteration. Subgroups of $\Out(F_r)$ that consist of only polynomially
growing outer automorphisms are known as \emph{Kolchin} subgroups and are
the $\Out(F_r)$ analog of unipotent subgroups of a linear group.

Clay and Pettet~\cite{clay-pettet} and Gultepe~\cite{gultepe} give sufficient
conditions for when a subgroup of $\Out(F_r)$ generated by two Dehn twists
contains an exponentially growing outer automorphism. This article gives an
algorithmic criterion that characterizes when a subgroup of $\Out(F_r)$
generated by two Dehn twists is Kolchin, in terms of a combinatorial
invariant of the gnerators known as the \emph{edge-twist} directed graph (see
Definition~\ref{def:et}).

\begin{theorem}[Main Theorem]
	Suppose $\sigma, \tau \in \Out(F_r)$ are Dehn twists. The subgroup
	$\langle \sigma,\tau\rangle$ is Kolchin if and only if the edge
	twist digraph of the defining graphs of groups is directed acyclic.
\end{theorem}

\section{Background}

\subsection{Graphs of groups and Dehn twists}

A \emph{graph} $\Gamma$ is a collection of vertices $V(\Gamma)$, edges
$E(\Gamma)$, initial and terminal vertex maps $o, t: E\to V$ and an involution
$\bar{\cdot}: E\to E$ satisfying $\bar{e}\neq e$ and $o(\bar{e}) = t(e)$. A
\emph{directed graph (digraph)} omits the involution. 

\begin{definition} 
	A \emph{graph of groups} is a pair $(G,\Gamma)$ where
	$\Gamma$ is a connected graph and $G$ is an assignment of groups to the
	vertices and edges of $\Gamma$ satisfying $G_e = G_{\bar{e}}$ and
	injections $\iota_e: G_e\to G_{t(e)}$. The assignment will often be
	suppressed and $\Gamma_v, \Gamma_e$ used instead.
\end{definition}

\begin{definition} 
	The \emph{fundamental groupoid} $\pi_1(\Gamma)$ of a graph
	of groups $\Gamma$ is the groupoid with vertex set $V(\Gamma)$
	generated by the path groupoid of $\Gamma$ and the groups $G_v$ subject
	to the following relations. We require that for each $v\in V(\Gamma)$
	the group $G_v$ is a subgroupoid based at $v$ and that the group and
	groupoid structures agree. Further, for all $e\in E(\Gamma)$ and $g\in
	G_e$ we have $\bar{e}\iota_{\bar{e}}(g)e = \iota_e(g)$.
\end{definition}

The \emph{fundamental group} of $\Gamma$ based at $v$, $\pi_1(\Gamma, v)$ is
the vertex subgroup of $\pi_1(\Gamma)$ based at $v$. It is standard that
changing the basepoint gives an isomorphic group~\cites{trees, higgins}. 

Let $(e_1,\ldots,e_n)$ be a possibly empty edge path in $\Gamma$ starting at
$v$ and $(g_0,\ldots, g_n)$ be a sequence of elements $g_i\in G_{t(e_i)}$ with
$g_0 \in G_v$. These data represent an arrow of $\pi_1(\Gamma)$ by the groupoid
product \[ g_0e_1g_1\cdots e_ng_n.\] A non-identity element of $\pi_1(\Gamma)$
expressed this way is \emph{reduced} if either $n = 0$ and $g_0\neq \id$, or $n
> 0$ ans for all $i$ such that $e_i = \bar{e}_{i+1}$, $g_i \not\in
\iota_{e_i}(G_{e_i})$. The edges appearing in a reduced arrow are uniquely
determined.
Further, if $t(e_n) = o(e_1)$ the arrow is cyclically reduced if either $e_n
\neq \bar{e}_0$ or $e_n = \bar{e}_0$ and $g_ng_0\not\in \iota_{e_n}(G_{e_n})$.
For an element $g\in \pi_1(\Gamma, v)$, the edges appearing in a cyclically
reduced arrow conjugate to $g$ in $\pi_1(\Gamma)$ is a conjugacy class
invariant.\footnote{These edges are covered by the axis of $g$ in the
Bass-Serre tree of $\Gamma$.}

\begin{definition} 
	Given a graph of groups $\Gamma$ a subset of edges
	$E'\subseteq E(\Gamma)$ and edge-group elements $\{z_e\}_{e\in E'}$
	satisfying $z_e \in Z(\Gamma_e)$ and $z_{\bar{e}} = z_e^{-1}$, the
	\emph{Dehn twist} of $\Gamma$ about $E'$ by $\{z_e\}$ is the
	fundamental groupoid
	automorphism $D_z$ given on the generators by
	\begin{align*}
		D_z(e) &= ez_e & e\in E' \\
		D_z(e) &= e & e\not\in E' \\
		D_z(g) &=g & g\in \Gamma_v \\
	\end{align*}
	As this groupoid automorphism preserves vertex subgroups it induces a
	well-defined outer automorphism class $D_z\in\Out(\pi_1(\Gamma, v))$,
	which we will also refer to as a Dehn twist.
\end{definition}

For a group $G$ we say $\sigma\in\Out(G)$ is a \emph{Dehn twist} if it can be realized
as a Dehn twist about some graph of groups $\Gamma$ with $\pi_1(\Gamma, v)
\cong G$.

Specializing to $\Out(F_r)$, when $\sigma\in\Out(F_r)$ is a Dehn twist there
are many graphs of groups $\Gamma$ with $\pi_1(\Gamma, v) \cong F_r$ that can
be used to realize $\sigma$. However, Cohen and Lustig~\cite{cohen-lustig}
define the notion of an \emph{efficient graph of groups} representative of a
Dehn twist and show that each Dehn twist in $\Out(F_r)$ has a unique efficient
representative. For a fixed $\sigma$ let $\mathcal{G}(\sigma)$ denote the graph
of groups of its efficient representative; edge groups of
$\mathcal{G}(\sigma)$ are infinite cyclic~\cite{cohen-lustig}. 

\begin{remark}\label{rem:twist-power}
	If $\sigma,\tau\in\Out(F_r)$ are Dehn twists with a common power, then
	$\mathcal{G}(\sigma) = \mathcal{G}(\tau)$.
\end{remark}

\subsection{Topological representatives and the Kolchin theorem}

Given a graph $\Gamma$ the \emph{topological realization} of $\Gamma$ is a
simplicial complex with zero-skeleton $V(\Gamma)$ and one-cells joining $o(e)$
and $t(e)$ for each edge in a set of $\bar{\cdot}$ orbit representatives. It
will not cause confusion to use $\Gamma$ for both a graph and its topological
representative. If $\gamma \subset \Gamma$ is a based loop, denote the
associated element of $\pi_1(\Gamma)$ by $\gamma^\ast$. Given $\sigma\in
\Out(F_r)$, a \emph{topological realization} is a homotopy equivalence
$\hat{\sigma}:\Gamma\to\Gamma$ so that $\hat{\sigma}_\ast : \pi_1(\Gamma,v)\to
\pi_1(\Gamma,\hat{\sigma}(v))$ is a representative of $\sigma$. A homotopy
equivalence $\hat{\sigma}:\Gamma\to\Gamma$ is \emph{filtered} if there is a
filtration $\emptyset = \Gamma_0\subsetneq \Gamma_1\subsetneq \cdots\subsetneq
\Gamma_k = \Gamma$ preserved by $\hat{\sigma}$.

\begin{definition}
	A filtered homotopy equivalence $\hat{\sigma}:\Gamma\to\Gamma$ is
	\emph{upper triangular} if
	\begin{enumerate}
		\item $\hat{\sigma}$ fixes the vertices of $\Gamma$,
		\item Each stratum of the filtration $\Gamma_i\setminus
			\Gamma_{i-1} = e_i$ is a single topological edge,
		\item Each edge $e_i$ has a preferred orientation and with this
			orientation there is a closed path $u_i\subseteq
			\Gamma_{i-1}$ based at $t(e_i)$ so that
			$\hat{\sigma}(e_i) = e_iu_i$.
	\end{enumerate}
\end{definition}

The path $u_i$ is called the suffix associated to $u_i$. A filtration assigns
each edge a height, the $i$ such that $e\in \Gamma_i\setminus\Gamma_{i-1}$, and
taking a maximum this definition extends to edge paths.

Every Dehn twist in $\Out(F_r)$ has an upper-triangular
representative~\cites{bfh-ii, cohen-lustig}. In a previous
paper~\cite{bering-lg} I describe how to construct $\mathcal{G}(\sigma)$ from
an upper-triangular representative, following a similar construction of
Bestvina, Feighn, and Handel~\cite{bfh-ii}. The following is an immediate
consequence of my construction.

\begin{lemma}\label{lem:filter}
	Suppose $\Gamma$ is a filtered graph and $\sigma\in\Out(F_r)$ is a Dehn
	twist that is upper triangular with respect to $\Gamma$. Then
	\begin{enumerate}
		\item there is a height function $ht: E(\mathcal{G}(\sigma))\to
			\mathbb{N}$ so that for any loop
			$\gamma\subseteq\Gamma_i$ the height of the edges in a
			cyclically reduced representative of the conjugacy
			class of $\gamma^\ast$ in
			$\pi_1(\mathcal{G}(\sigma))$ is at most $i$,
		\item For each edge $e\in E(\mathcal{G}(\sigma))$, $ht(e) =
			ht(\bar{e})$, and the edge
			group $\mathcal{G}(\sigma)_e$ is a conjugate of a
			maximal cyclic subgroup of $F_r$ that contains
			$u_i^\ast$ for some suffix $u_i$, and $ht(e) >
			\min_{\gamma \sim u_i} \{ht_\Gamma(\gamma)\}$, where
			$\gamma$ ranges over loops freely homotopic to $u_i$.
	\end{enumerate}
\end{lemma}

Bestvina, Feighn, and Handel proved an $\Out(F_r)$ analog of the classical
Kolchin theorem for $\Out(F_r)$, which provides simultaneous upper-triangular
representatives for Kolchin-type subgroups of $\Out(F_r)$.

\begin{theorem}[\cite{bfh-ii}]\label{thm:kolchin}
	Suppose $H\leq\Out(F_r)$ is a Kolchin subgroup. Then there is a finite
	index subgroup $H' \leq H$ and a filtered graph $\Gamma$ so that each
	$\sigma \in H'$ is upper triangular with respect to $\Gamma$.
\end{theorem}

\subsection{Twists and polynomial growth}

In a previous paper~\cite{bering-lg} I introduced the edge-twist digraph of two
Dehn twists and used it to provide a sufficient condition for a subgroup of
$\Out(F_r)$ generated by two Dehn twists to be Kolchin.

\begin{definition}[\cite{bering-lg}]\label{def:et}
	The edge-twist digraph $\mathcal{ET}(A,B)$ of two graphs of groups $A$,
	$B$ with isomorphic fundamental groups and infinite cyclic edge
	stabilizers is the digraph with vertex set
	\[ V(\mathcal{ET}(A,B)) = \{(e, \bar{e}) | e\in E(A)\}\cup
	\{(f,\bar{f})|f\in E(B)\} \]
	directed edges $((e,\bar{e}), (f,\bar{f}))$ $e\in E(A), f\in E(B)$ when a generator of $A_e$
	uses $f$ or $\bar{f}$ in its cyclically reduced representation in
	$\pi_1(B)$, and directed edges $((f,\bar{f}),(e,\bar{e}))$, $f\in E(B),
	e\in E(A)$ when a generator of $B_f$ uses $e$ or $\bar{e}$ in a
	cyclically reduced representation in $\pi_1(A)$.
\end{definition}

\begin{remark}
	This is well-defined, using an edge is a conjugacy invariant, and using
	an edge or its reverse is preserved under taking inverses.
\end{remark}

\begin{lemma}[\cite{bering-lg}]\label{lem:pg}
	If $\sigma, \tau\in\Out(F_r)$ are Dehn twists and
	$\mathcal{ET}(\mathcal{G}(\sigma),\mathcal{G}(\tau))$ is directed acyclic, then 
	$\langle \sigma,\tau\rangle$ is Kolchin.
\end{lemma}

\section{Proof of the main theorem}

\begin{proof}
	It suffices to prove the converse to Lemma \ref{lem:pg}. Suppose 
	$\langle \sigma, \tau\rangle$ is Kolchin.
	By Theorem \ref{thm:kolchin}, there is a finite
	index subgroup $H\leq \langle \sigma, \tau\rangle$
	where every element of $H$ is upper triangular with respect
	respect to a fixed filtered graph $\Gamma$. Since $H$ is finite index, there are
	powers $m, n$ so that $\sigma^m, \tau^n \in H$, so that $\sigma^m$
	and $\tau^n$ are upper triangular with respect to $\Gamma$. 

	By Lemma \ref{lem:filter}, the filtration of $\Gamma$ induces height
	functions on $E(\mathcal{G}(\sigma^m))$ and $E(\mathcal{G}(\tau^n))$,
	combining these gives a height function on the vertices of
	$\mathcal{ET}(\mathcal{G}(\sigma^m),\mathcal{G}(\tau^n))$. Every
	directed edge $((e,\bar{e}), (f,\bar{f}))$ satisfies $ht(e) > ht(f)$.
	Indeed, suppose $(e,\bar{e}) \in E(\mathcal{G}(\sigma^m))$. Let
	$[g] \subset F_r$ be the conjugacy class of a generator of
	$\mathcal{G}(\sigma^m)_e$. By
	Lemma \ref{lem:filter} (ii), there is a representative $g\in [g]$ such that
	$g^k = u_i^\ast$ for some $\sigma$-suffix $u_i$. Take a minimum
	height loop $\gamma$ representing $[g]$, so that $\gamma \subseteq
	\Gamma_{ht(\gamma)}$. Again by Lemma \ref{lem:filter} (ii), $ht(e) >
	ht(\gamma)$. Finally, by Lemma \ref{lem:filter} (i), each edge $f$ in a
	cyclically reduced $\pi_1(\mathcal{G}(\tau^n))$ representative of
	$[g]$ satisfies $ht(f) \leq ht(\gamma)$. Thus $ht(e) > ht(f)$ for each
	directed edge with origin $(e,\bar{e})$, as required. The argument for
	generators of the edge groups of
	$\mathcal{G}(\tau^n)$ is symmetric. Therefore any directed path
	in $\mathcal{ET}(\mathcal{G}(\sigma^m),\mathcal{G}(\tau^n))$ has
	monotone decreasing vertex height, which implies that $\mathcal{ET}$
	is directed acyclic. To conclude, observe that by
	Remark~\ref{rem:twist-power}
	$\mathcal{ET}(\mathcal{G}(\sigma^m),\mathcal{G}(\tau^n)) =
	\mathcal{ET}(\mathcal{G}(\sigma), \mathcal{G}(\tau))$.
\end{proof}

Cohen and Lustig~\cite{cohen-lustig-ii} give an algorithm to find efficient
representatives. Computing the edge-twist digraph and testing if it is acyclic
are straightforward computations, so the criterion in the main theorem is
algorithmic.

\section*{Acknowledgements}

I thank Mladen Bestvina for a helpful conversation that led to this note.

\end{document}